\documentclass[reqno,12pt,b5paper]{amsart}
\usepackage[lmargin=2.5cm,rmargin=2.5cm,bmargin=2cm,tmargin=2cm]{geometry}
\usepackage{amssymb}
\usepackage{esint}
\usepackage{mathtools}
\usepackage[bookmarksopen,bookmarksdepth=3]{hyperref}
\usepackage[shortlabels]{enumitem}
\usepackage[final]{pdfpages}
\usepackage[finnish,english]{babel}
\usepackage[utf8]{inputenc} 

\theoremstyle{plain}
\newtheorem{theorem}{Theorem}[section]

\newtheorem{corollary}[theorem]{Corollary}

\theoremstyle{definition}

\theoremstyle{remark}
\newtheorem{remark}[theorem]{Remark}

\newcommand{\vnabla}{\overset{v}{\nabla}}
\newcommand{\hnabla}{\overset{h}{\nabla}}
\newcommand{\T}{\mathbb{T}}
\newcommand{\R}{\mathbb{R}}
\newcommand{\Z}{\mathbb{Z}}
\newcommand{\C}{\mathbb{C}}
\newcommand{\N}{\mathbb{N}}
\newcommand{\Q}{\mathbb{Q}}
\newcommand{\Hset}{\mathbb{H}}

\newcommand{\Gr}{\mathbf{Gr}}

\newcommand{\vev}[1]{\left\langle#1\right\rangle}

\newcommand{\eps}{\varepsilon}
\renewcommand{\phi}{\varphi}
\newcommand{\abs}[1]{\left| #1 \right|}
\newcommand{\norm}[1]{\Vert #1 \Vert}
\newcommand{\aabs}[1]{\left\| #1 \right\|}

\DeclareMathOperator*{\argmin}{arg\,min}

\begin{document}
\pagenumbering{roman}
\setcounter{page}{1}
\thispagestyle{plain}
\begin{center}
	~\\
	\vspace{2.5cm}
	\textsc{\Large
		{Geodesic tomography problems on Riemannian manifolds}\\
		\vspace{1.5cm}
		\large
		Jesse Railo
		\normalsize
	}
\end{center}
\newpage

\title[Geodesic tomography problems on Riemannian manifolds]{}

\author[Jesse Railo]
\maketitle
\newpage
\thispagestyle{plain}
Editors 

\textit{Mikko Salo}

Department of Mathematics and Statistics

University of Jyväskylä

\newpage
\thispagestyle{plain}

\section*{Foreword} 
I wish to thank my advisor, Mikko Salo, for his support and help during my PhD studies at the University of Jyväskylä. He has been the best teacher and academic role model that I could have hoped for. I owe him a great debt of gratitude. I thank the Department of Mathematics and Statistics for giving me a friendly working environment and support in 2015--2019.

I wish to thank my collaborators, Joonas Ilmavirta, Olli Koskela, and Jere Lehtonen, for many productive and instructive discussions. I have learned many good ways of thinking and working from you. I also thank any other colleagues I have interacted with and whom are not mentioned here by name.

I wish to thank François Monard who has agreed to be the opponent at the public examination of my dissertation. I wish to thank Todd Quinto and Hanming Zhou for their preliminary examinations of my dissertation, which I have already received when writing this.

Finally, I thank my wife, Heli, for supporting and loving me during these years. She has been very compassionate for my, sometimes, comical working hours and habit to work at home. I thank my family and friends for offering me many great opportunities to take a break from mathematics to do something completely different and fun. \linebreak \linebreak

Jyväskylä, November 11, 2019

Department of Mathematics and Statistics

University of Jyväskylä

\textit{Jesse Railo}

\newpage
\thispagestyle{plain}

\section*{List of included articles} 
This dissertation consists of an introductory part and the following four articles: \\

\begin{enumerate}[label=(\Alph*)]
\item \label{paper1} Jere Lehtonen, Jesse Railo and Mikko Salo. \textit{Tensor tomography on Cartan-Hadamard manifolds.} Inverse Problems 34 (2018), special issue: 100 years of the Radon transform, no. 4, 044004. \linebreak
\item \label{paper2} Joonas Ilmavirta and Jesse Railo. \textit{Geodesic ray transform with matrix weights for piecewise constant functions.} Preprint (2019), arXiv:1901.03525.\linebreak
\item \label{paper3} Joonas Ilmavirta, Olli Koskela and Jesse Railo. \textit{Torus computed tomography.}\,\,Preprint\,(2019),\,\,arXiv:1906.05046.\linebreak
\item \label{paper4} Jesse Railo. \textit{Fourier analysis of periodic Radon transforms.} Preprint (2019), arXiv:1909.00495.
\end{enumerate}
\hfill 

The author of this dissertation has actively taken part in the research of the joint articles \ref{paper1}, \ref{paper2} and \ref{paper3}.

\newpage
\thispagestyle{plain}

\section*{Tiivistelmä}

\begin{otherlanguage}{finnish}Väitöskirjassa tutkitaan integraaligeometriaan liittyviä inversio-ongel-mia. Geodeettinen sädemuunnos on operaattori, joka laskee funktion polkuintegraalin geodeesia pitkin. Väitöskirjassa määritetään monia ehtoja, joilla tällainen tieto määrää funktion yksikäsitteisesti ja vakaasti. Lisäksi osana väitöskirjan työtä on toteutettu numeerinen malli, jota voidaan käyttää tietokonetomografiassa.

Väitöskirjan johdannossa esitetään inversio-ongelmien peruskäsitteitä ja tietokonetomografiaan läheisesti liittyviä matemaattisia malleja. Johdannon pääpaino on integraaligeometriaan liittyvien mallien määrittelys-sä, tutkimusaiheen kirjallisuuskatsauksessa ja väitöskirjan tutkimustulosten esittelyssä. Lisäksi annetaan lista integraaligeometrian tärkeistä avoimista matemaattisista ongelmista.

Väitöskirjan ensimmäisessä artikkelissa osoitetaan, että symmetrinen solenoidaalinen tensorikenttä voidaan määrätä yksikäsitteisesti sen geodeettisesta sädemuunnoksesta Cartan-Hadamard monistolla, kun tietyt geometriasta riippuvat vähenemisehdot täyttyvät. Tutkittu integraalimuunnos esiintyy sirontaan liittyvissä käänteisongelmissa kvanttifysiikassa ja yleisessä suhteellisuusteoriassa.

Väitöskirjan toisessa artikkelissa näytetään, että paloittain vakio vektoriarvoinen funktio voidaan määrittää yksikäsitteisesti sen matriisipainotetusta geodeettisesta sädemuunnoksesta reunallisella Riemannin monistolla, jos geometria sallii aidosti konveksin funktion olemassaolon ja epäsingulaarinen matriisipaino riippuu jatkuvasti sen sijainnista moniston yksikköpallokimpulla. Tällaista integraalimuunnosta voidaan käyttää mallintamaan attenuoitua sädemuunnosta sekä inversio-ongelmia konnektiolle ja Higgsin kentälle.

Väitöskirjan kolmannessa ja neljännessä artikkelissa tutkitaan geodeettista sädemuunnosta suljettujen geodeesien yli toruksella, kun funktioiden säännöllisyys on alhainen. Neljännessä artikkelissa tarkastellaan lisäksi tällaisen muunnoksen yleistystä, kun funktion integraalit tunnetaan isometrisesti upotettujen alempiasteisten toruksien yli. Artikkeleissa todistetaan uusia rekonstruktiokaavoja, regularisointistrategioita ja vakausestimaatteja tällaisille integraalimuunnoksille. Saaduilla tutkimustuloksilla on sovelluskohteita erilaisissa laskennallisissa tomografiamenetelmissä.
\end{otherlanguage}
\newpage
\thispagestyle{plain}
\section*{Abstract} 
This dissertation is concerned with integral geometric inverse problems. The geodesic ray transform is an operator that encodes the line integrals of a function along geodesics. The dissertation establishes many conditions when such information determines a function uniquely and stably. A new numerical model for computed tomography imaging is created as a part of the dissertation.

The introduction of the dissertation contains an introduction to inverse problems and mathematical models associcated to computed tomography. The main focus is in definitions of integral geometry problems, survey of the related literature, and introducing the main results of the dissertation. A list of important open problems in integral geometry is given.

In the first article of the dissertation, it is shown that a symmetric solenoidal tensor field can be determined uniquely from its geodesic ray transform on Cartan-Hadamard manifolds, when certain geometric decay conditions are satisfied. The studied integral transforms appear in inverse scattering theory in quantum physics and general relativity.

In the second article of the dissertation, it is shown that a piecewise constant vector-valued function can be determined uniquely from its geodesic ray transform with a continuous and non-singular matrix weight on Riemannian manifolds that admit a strictly convex function and have a strictly convex boundary. These integral transforms can be used to model attenuated ray transforms and inverse problems for connections and Higgs fields.

The third and fourth articles of the dissertation study the geodesic ray transform over closed geodesics on flat tori when the functions have low regularity assumptions. The fourth article studies a generalization of the geodesic ray transform when the integrals of a function are known over lower dimensional isometrically embedded flat tori. New inversion formulas, regularization strategies and stability estimates are proved in the articles. The new results have applications in different computational tomography methods.
\newpage

\pagenumbering{arabic}
\setcounter{tocdepth}{1}
\tableofcontents

\section{Introduction}

One of the most fundamental inverse problems asks if an unknown function is determined uniquely from the knowledge of the values of its line integrals over all possible lines in Euclidean space. This is in fact the mathematical model used for X-ray computed tomography (CT). This can be viewed as an integral transform acting on functions. Its many possible generalizations model other tomographic methods such as computerized axial tomography (CAT), positron-emission tomography (PET) and single-photon emission tomography (SPECT). It also has close connection to other inverse problems and applications such as seismic imaging, electrical impedance tomography, polarization tomography, quantum state tomography, inverse spectral problems and inverse scattering problems. This thesis studies generalizations of X-ray tomography on Riemannian manifolds.

This introductory part of the thesis is organized as follows. We discuss inverse problems and X-ray computed tomography in general in sections \ref{sec:invprob} and \ref{sec:xrays} respectively. We shortly describe the articles \ref{paper1}--\ref{paper4} in section \ref{sec:summary}. Preliminaries on Riemannian manifolds are given in section \ref{sec:riemannian}. We define different geodesic tomography models and corresponding inverse problems in section \ref{sec:results}. We also survey related solved and unsolved problems in section \ref{sec:results}. We introduce the main results of this thesis in section \ref{sec:main}. The results are proved in the included articles \ref{paper1}--\ref{paper4}.

\subsection{Inverse problems}\label{sec:invprob}

Inverse problems is a field of mathematics where one typically measures data outside or on the boundary of an object and wants to recover knowledge of its internal structure. Such mathematical problems occur often in medical, engineering and physical applications. In some inverse problems, measurements are done very far from an object. Such problems can be naturally studied using noncompact spaces in mathematical models.

Typical mathematical questions that one studies in inverse problems include:

\begin{enumerate}[i)]
\item \textit{(Forward problem)} What is a good mathematical model that captures the physical phenomenon which relates measurement data to physical parameters of an unknown object? Does the mathematical model define data uniquely? \label{eka}
\item \textit{(Uniqueness)} Do measurements determine the unknown physical parameters uniquely? If not, can non-uniqueness be characterized? \label{toka}
\item \textit{(Reconstruction)} How can the unknown physical parameters be computed from measurement data? \label{kolm}
\item \textit{(Stability)} Do the unknown physical parameters depend continuously on measurement data? Does there exists a quantitative stability estimate? \label{nelj}
\item \textit{(Simulations and regularization)} How can reconstruction methods be implemented into numerical algorithms? How to overcome instability caused by ill-posedness and measurement noise, finiteness of measurements, and numerical approximations? \label{vika}
\end{enumerate}

The questions \ref{eka}--\ref{toka} are encountered in \ref{paper1}, the question \ref{toka} in \ref{paper2}, and the questions \ref{kolm}--\ref{vika} in \ref{paper3} and \ref{paper4}. The textbooks \cite{KS05, K11, MS12} and the survey \cite{U14} can be used to find more details and references on inverse problems in general.

\subsection{X-ray tomography and its generalizations}\label{sec:xrays} Let $f$ be a function $\R^n \to \R$. One defines the \textit{X-ray transform} in $\R^n$ as \begin{equation}\label{eq:xray}\mathcal{R}f(x,v) = \int_{\R} f(x+tv)dt\end{equation} where $(x,v) \in \R^n \times S^{n-1}$ whenever the integral is well-defined and finite. This is the standard mathematical model for X-ray tomography measurements, and it is also known as the \textit{Radon transform} if $n = 2$. In higher dimensions, the X-ray transform and the Radon transform are different operators \cite{H99}. The corresponding uniqueness problem asks if $\mathcal{R}f = \mathcal{R}g$ implies that $f = g$. The other questions \ref{eka}--\ref{vika} of section \ref{sec:invprob} could be asked as well.

The inverse problems associated to the X-ray transform were first studied by Johann Radon in 1917 \cite{R17}. Fritz John characterized the range of the X-ray transform in $\R^3$ in terms of ultrahyperbolic equations (called \textit{John's equations}) in 1938 \cite{J38}. Later, the mathematical problem was restudied independently by Allan Cormack in 60s \cite{C63,C64}. Godfrey Hounsfield studied practical CT imaging a few years later. For their seminal works on CT imaging, Cormack and Hounsfield won the 1979 Nobel Prize in Physiology or Medicine. The monographs \cite{N01, KS01, H99} and the surveys \cite{Q06,Hel13,KQ15} are recommended references on the mathematics of the X-ray and Radon transforms.

The X-ray transform can be generalized many ways:

\begin{enumerate}[i)]
\item Instead of integrating over straight lines, suppose one knows integrals of $f$ over other families of curves. For example, data could be measured over geodesics of a Riemannian manifold. \label{item:geod}
\item Instead of integrating against the measure $dt$, suppose one knows integrals of $f$ against the weighted measure $w(x,v)dt$ where $w(x,v) > 0$ is a continuous function on $\R^n \times S^{n-1}$. \label{item:weighted}
\item Instead of integrating over straight lines, suppose one knows integrals of $f$ over other families of sets. For example, data could be measured over hyperplanes.\label{item:radon}
\item Instead of integrating a function, suppose one knows integrals of a tensor field so that the value of $f$ depends also on the direction of an X-ray, not only on a point in $\R^n$.\label{item:tensor}
\item Some combination of the above cases.\label{item:mix}
\end{enumerate}

The case \ref{item:geod} corresponds to the \textit{geodesic X-ray transform}, \ref{item:weighted} to the \textit{X-ray transform with weights}, \ref{item:radon} to the \textit{Radon transform}, and \ref{item:tensor} to \textit{tensor tomography}. These different generalizations of the X-ray transform are studied in this thesis. One of the fundamental properties is that all of these integral transforms are linear. This reduces the uniqueness problem to studying kernels of the transforms. 

The field of inverse problems that studies these integral transforms, among other problems of similar nature, is often called \textit{integral geometry}. For example, the boundary rigidity problem asks if the knowledge of distances between any two boundary points determines the geometrical shape of a compact connected object with boundary uniquely (see section \ref{sec:openproblems} for a rigor formulation). This is an example of a nonlinear integral geometry problem. We give a more detailed introduction to integral geometry problems in section \ref{sec:results}. More references and recent developments in integral geometry can be found from the textbook \cite{S94} and the surveys \cite{PSU14,IM18}.

\subsection{On the articles in this thesis} \label{sec:summary}

The first article \ref{paper1} with Lehtonen and Salo considers tensor tomography on Cartan-Hadamard manifolds. Tensors can be used for modeling physical parameters that have spatial and directional dependence. In this work, we characterize the kernel of the geodesic ray transform for symmetric tensor fields of any order under sufficient decay conditions. This generalizes injectivity results of the geodesic ray transform from compact manifolds with boundary to noncompact manifolds.

The second article \ref{paper2} with Ilmavirta considers the geodesic ray transform with matrix weights on manifolds that admit a strictly convex function. In this work, we restrict our study to the class of piecewise constant vector-valued functions. We show injectivity of this transform under the assumption that the weight is continuous and invertible at any point. This assumption on weights is very mild, and counterexamples for injectivity on smooth functions exist even in Euclidean case. The geometric assumption is equivalent to a manifold being nontrapping in dimension two. Injectivity of the geodesic ray transform (without a weight) for smooth functions on nontrapping manifolds is one of the most important unanswered geometric inverse problem at the moment.

The third article \ref{paper3} with Ilmavirta and Koskela studies the geodesic X-ray transform over periodic geodesics on the flat $2$-torus. In this work, reconstruction methods, including regularization and numerical implementations, drive theoretical considerations. We prove new reconstruction formulas for integrable functions, solve a minimization problem associated to Tikhonov regularization in Sobolev spaces, and prove that the unique minimizer provides a regularization strategy. We have also computed and analyzed the adjoint and the normal operators. Regularization of reconstructions is important since measurement noise is amplified in practice due to ill-posedness of the problem. Another reason for regularization is that one can collect only finitely many measurements in practice. We created Matlab codes, performed numerical tests and demonstrated how the developed methods can be applied in practical CT imaging.

The fourth article \ref{paper4} studies the $d$-plane Radon transforms on the flat $n$-tori $\T^n$. The main results in \ref{paper4} extend theorems in \ref{paper3} to higher dimensions. In addition, new stability estimates in Bessel potential norms and inversion formulas for periodic distributions are proved. It is shown that the $d$-plane Radon transforms maps the Bessel potential spaces continuously into the weighted Bessel potential spaces on $\T^n \times \Gr(d,n)$ where $\Gr(d,n)$ is the collection of $d$-dimensional subspaces of $\Q^n$. The use and analysis of such structures is the main methodological advance compared to \ref{paper3}. Quite surprisingly, one of the inversion formulas in \ref{paper4} implies that a compactly supported function on the plane with zero average is a sum of its X-ray data.

\section{Preliminaries on Riemannian manifolds}\label{sec:riemannian}

Let $(M,g)$ denote a Riemannian manifold with or without boundary. We assume always that $M$ is complete and $\dim(M) \geq 2$. We define the following notations: \begin{itemize}
\item The unit tangent bundle is denoted by \begin{equation}SM = \{\,(x,v) \in TM \,;\, \abs{v}_g = 1\,\}.\end{equation}
\item If $(x,v) \in SM$, then $\gamma_{x,v}$ denotes the unique unit-speed geodesic such that $\gamma(0) = x$ and $\dot{\gamma}(0) = v$. The set of maximal unit-speed geodesics of $M$ is denoted by $\Gamma$.
\item We denote the boundary of $M$ by $\partial M$ and by $\nu(x)$ the inward pointing unit normal of $\partial M$ at $x \in \partial M$.
\item We say that $M$ has a \textit{strictly convex boundary} if the second fundamental form of $\partial M$ is positive definite or, equivalently, principal curvatures of $\partial M$ are positive.
\item We denote the covariant derivative by $\nabla$ and the Riemannian curvature tensor by $R$.
\item We denote the sectional curvature of a two-plane $\Pi \subset T_xM$ by $K_x(\Pi)$ and $\mathcal{K}(x) = \sup \{\,\abs{K_x(\Pi)}\,;\, \Pi \subset T_xM \text{ is a two-plane}\,\}.$ 
\item We write $K \leq 0$ if $K_x(\Pi) \leq 0$ for any $x \in M$ and any two-plane $\Pi \subset T_xM$. In this case, we say that $M$ has \textit{non-positive (sectional) curvature}.
\end{itemize}

\subsection{Definitions related compact Riemannian manifolds with boundary}

Let $(M,g)$ be a compact Riemannian manifold with a strictly convex boundary. We define some useful geometric terminology in this sections. In the following sections, we give results on geodesic ray transforms using these different geometric definitions.

We say that $M$ is \textit{simple} if the exponential map $\exp_p: T_pM \to M$ is a diffeomorphism from its maximal domain for any $p \in M$. This, in particular, implies that there are no conjugate points, any two points are connected by a unique geodesic, and $M$ is diffeomorphic to the Euclidean ball of dimension $\dim(M)$ \cite{Pet06}. We say that $M$ is \textit{nontrapping} if $\gamma_{x,v}(t)$ meets $\partial M$ in finite time for any $(x,v) \in SM$. In particular, simple manifolds are nontrapping.

We say that $f: M \to \R$ is a \textit{strictly convex function} if $f \in C^\infty(M)$ so that $\text{Hess}_x(f)$ is positive definite for any $x \in M$ or, equivalently, $(f\circ \gamma)''(t) > 0$ for every geodesic $\gamma \in \Gamma$. A manifold $M$ satisfies the \textit{foliation condition} if there exists a strictly convex function \cite{UV16,PSUZ16}. 

\begin{remark} The level sets of a strictly convex function are strictly convex hypersurfaces besides the special case of the minimum whose level set is a single point \cite{PSUZ16}. The corresponding level sets form layers that foliate the whole manifold. The tangential geodesics of a strictly convex hypersurface do not locally travel inside the hypersurface. Using the foliation condition, this type of behavior can be made global. In turn, this allows one to use a layer stripping argument for proving injectivity of the geodesic ray transform if local injectivity can be shown \cite{UV16}.
\end{remark}

The \textit{trapped set of $M$}, denoted by $K \subset SM$, consists of points $(x,v) \in SM$ such that $\gamma_{x,v}(t)$ does not meet the boundary $\partial M$ for any $t \in \R$. In particular, if $M$ is nontrapping (in the sense of above), then $K = \emptyset$. The trapped set is said to be \textit{hyperbolic} if there is a certain orthogonal splitting to geodesic, stable and unstable parts of $T_{(x,v)}(SM)$ for any $(x,v) \in K$. For exact definitions, see \cite{Pat99,G17}.

Let $\beta \geq 0$. We say that $J$ is \textit{$\beta$-Jacobi field} along $\gamma \in \Gamma$ if it satisfies
\begin{equation} D_t^2 J(t) + \beta R(J(t),\dot{\gamma}(t))\dot{\gamma}(t) = 0.
\end{equation} We say that two distinct points along $\gamma$ are \textit{$\beta$-conjugate} if there exists a non-trivial $\beta$-Jacobi field which vanishes at the points. The \textit{$\beta$-terminator value} $\beta_{\text{Ter}}$ is the supremum of the numbers $\beta$ so that $M$ is free of $\beta$-conjugate points. In particularly $\beta_{Ter} = \infty$ if and only if $K \leq 0$, and $M$ has no conjugate points if and only if $\beta_{Ter} \geq 1$. For more details see \cite{PSU15}.

\begin{remark} We do not study manifolds that have trapped geodesics in this thesis, but this condition is included as we will give references to other works on the geodesic ray transform where a hyperbolic trapped set is a part of the geometrical assumptions. Our reason for introducing $\beta$-Jacobi fields here is similar and they are not applied in this thesis.
\end{remark}

\subsection{Cartan-Hadamard manifolds}

We say that a Riemannian manifold $(M,g)$ without boundary is a \textit{Cartan-Hadamard manifold} if $(M,g)$ is complete, simply connected and $K \leq 0$. The classical Cartan-Hadamard theorem states that $\exp_p: T_pM \to M$ is a diffeomorphism for any $p \in M$ (see e.g. \cite[Chapter 6]{Pet06} or \cite[Chapter 11]{Lee97}). In particular, $M$ with $\dim(M) = n$ is diffeomorphic to Euclidean space $\R^n$. This implies that Cartan-Hadamard manifolds are noncompact.

The model spaces of Cartan-Hadamard manifolds are the hyperbolic space $\Hset^n$ ($K \equiv -1$) and Euclidean space $\R^n$ ($K \equiv 0$). Many other examples can be constructed using warped products with radial metrics \cite{BO69,KW74,GW79,Pet06}. A discussion on such constructions, related to the theorems of the article \ref{paper1}, is given in [\ref{paper1}, Section 2].

\section{Geodesic tomography problems}
\label{sec:results}

\subsection{Geodesic tensor tomography}

We denote by $C^1(T^mM)$ the set of $C^1$-smooth covariant $m$-tensor fields of $M$ and by $C^1(S^mM) \subset C^1(T^mM)$ the set of symmetric covariant $m$-tensor fields. Each $f \in C^1(T^mM)$ can be written in local coordinates as
\begin{equation}f = f_{j_1\cdots j_m}(x) dx^{j_1} \otimes \cdots \otimes dx^{j_m}\end{equation}
using the Einstein summation convention. Let $\Pi_M$ denote the permutation group of $\{1,\dots,m\}$. Tensors in $f \in C^1(S^mM)$ are symmetric in the sense that
\begin{equation}f = f_{j_{\sigma(1)} \cdots j_{\sigma(m)}}(x) dx^{j_1} \otimes \cdots \otimes dx^{j_m}\end{equation} for any $\sigma \in \Pi_m$.

If every maximal geodesic of $M$ has finite length, then one defines the \textit{geodesic ray transform of symmetric $m$-tensor fields} by the formula
\begin{equation}\label{eq:Imdef}I_mf(\gamma) = \int_{\gamma} \lambda_m f(\gamma(t),\dot{\gamma}(t))dt\end{equation}
where $\gamma \in \Gamma$ and $\lambda_m f(x,v) = f_{j_1\cdots j_m}(x)v^{j_1}\cdots v^{j_m}$ is a mapping $SM \to \R$. In fact, $\lambda_m$ maps $C^1(S^mM) \to C^1(SM)$ so that the spherical harmonics decomposition with respect to $v$ of $\lambda_m f$ is of degree $m$. A more detailed exposition of symmetric tensors and $\lambda_m$ are given in \cite{S94, DS11}. There is also a brief discussion in [\ref{paper1}, Section 3.3].

In general, the geodesic ray transform $I$ can be straightforwardly defined for every $f \in C(SM)$ if every maximal geodesic of $M$ has finite length. However, this transform always has a non-trivial kernel on manifolds with boundary, even in the case of symmetric $m$-tensor fields with $m \geq 1$, as we will explain later. This motivates to study functions of $C(SM)$ that have a special form. 

In the article \ref{paper1}, we study the kernel of the geodesic X-ray transform on Cartan-Hadamard manifolds for functions that arise from symmetric tensors. In this case, any maximal geodesic of $M$ has infinite length. Therefore, the integrals (\ref{eq:Imdef}) are finite only if the tensors decay sufficiently fast along every geodesic.

\subsubsection{On the kernel of $I_m$ and solenoidal injectivity} We define the s\textit{ymmetrization} of a tensor $\sigma_m: T^mM \to S^mM$ by
\begin{equation}\sigma_m(f) = \frac{1}{m!}\sum_{\sigma \in \Pi_m}f_{j_{\sigma(1)}\cdots j_{\sigma(m)}}(x) dx^{j_1} \otimes \cdots \otimes dx^{j_m}.\end{equation}


Let $\phi_t(x,v) = (\gamma_{x,v}(t),\dot{\gamma}_{x,v}(t))$ be the \textit{geodesic flow} on $SM$. One defines the \textit{geodesic vector field} $X$ for functions in $C^1(SM)$ as
\begin{equation}Xf(x,v) := \frac{d}{dt} f(\phi_t(x,v))|_{t=0}.\end{equation}

Suppose now that $M$ is a nontrapping Riemannian manifold. One can show that $IXf = 0$ for any $f \in C^1(SM)$ with $f|_{\partial M} = 0$ by the fundamental theorem of calculus. Another calculation shows that \begin{equation}X(\lambda_m f) = \lambda_m(\sigma_m \nabla f)\end{equation} for any $f \in C^1(S^mM)$. Therefore, if $m \geq 1$, the kernel of $I_m$ contains all symmetric $m$-tensors of the form $\sigma_m \nabla f$ where $f \in C^1(S^{m-1}M)$ and $f|_{\partial M}= 0$. We say that $f$ is a \textit{potential} of the tensor $\sigma_m \nabla f$.

We identify the space $C^1(S^{-1}M)$ as the space of the zero function. We say that $I_m$ is \textit{$s$-injective} if the kernel of $I_m$ contains only tensors that arise from a potential described above. This implies that the solenoidal part of a symmetric tensor can be uniquely determined from its geodesic ray transform (see \cite{S94} for details about the Helmholtz decomposition of symmetric tensors). We list next some known injectivity results for smooth tensor fields on compact Riemannian manifolds with a strictly convex boundary:

\begin{itemize}
\item If $M$ is a simple manifold, then $I_m$ is s-injective for $m = 0,1$ \cite{M77,AR97}.
\item If $M$ is a simple manifold whose metric is from a generic class (including real analytic metrics), then $I_2$ is s-injective \cite{SU05}.
\item If $M$ is a simple manifold of $\dim(M) = 2$, then $I_m$ is s-injective for every $m \geq 0$ \cite{PSU13}.
\item If $M$ is a nontrapping manifold of $\dim(M) = 2$, $I_m$ is s-injective for $m = 0,1$ and the adjoint of $I_0$ is surjective, then $I_m$ is s-injective for every $m \geq 0$ \cite{PSU13}.
\item If $M$ is a simple manifold with $n = \dim(M) \geq 2$ and $\beta_{Ter} \geq \frac{m(m+n-1)}{2m+n-2}$, then $I_m$ is s-injective \cite{PSU15}.
\item If $M$ is a nontrapping manifold of $\dim(M) \geq 3$ with a strictly convex foliation, then $I_m$ is s-injective for $m = 0,1,2,4$ \cite{UV16,SUV17,HUZ18}.
\item If $M$ is a compact Riemannian manifold with no conjugate points and hyperbolic trapped set, then $I_m$ is s-injective for $m = 0,1$. If moreover $K \leq 0$, then $I_m$ is s-injective for every $m \geq 0$ \cite{G17}.
\item If $M$ is a compact Riemannian manifold of $\dim(M) = 2$ with no conjugate points and hyperbolic trapped set, then $I_m$ is s-injective for every $m \geq 0$ \cite{L19}.
\item If $M$ is a simple manifold with real analytic metric, then $I_m, m \in \N$, admit a certain local support theorem \cite{KS09,AM19}. (These results are partly contained in the results of \cite{UV16,SUV17,HUZ18}.)
\end{itemize}

We state some of the related open problems in section \ref{sec:openproblems}.

\subsection{Geodesic ray transform with matrix weights}

Suppose that $W: SM \to \C^{m\times m}$ is continuous and $W(x,v): \C^m \to \C^m$ is injective for any $(x,v) \in SM$. Let $f: SM \to \C^m$ be a continuous function. One can then define \textit{the geodesic ray transform with the weight $W$} as
\begin{equation}I_Wf(x,v) := \int_{a_{x,v}}^{b_{x,v}} W(\gamma_{x,v}(t),\dot{\gamma}_{x,v}(t))f(\gamma_{x,v}(t),\dot{\gamma}_{x,v}(t))dt\end{equation} where $[a_{x,v},b_{x,v}]$ is the maximal domain of $\gamma_{x,v} \in \Gamma$ (possibly infinite).

The corresponding uniqueness problem asks if the knowledge of $I_Wf$ and $W$ determine $f$ uniquely. There exist counterexamples and positive results to the uniqueness problem. Clearly, if $W$ does not depend on the coordinate $v$, then injectivity of $I_W$ is equivalent to injectivity of $I$ without a weight (i.e. $W \equiv 1$).

An important special case of the geodesic ray transforms with weights is the \textit{attenuated geodesic ray transform}. The attenuated geodesic ray transforms is studied very recently for example in \cite{SU11,PSU12,AMU18,HMS18,MNP19a,MNP19b,BM19}. In the simplest model for the attenuated ray transform (with $m = 1$), the weight has a special form
\begin{equation}w_a(x,v) = \exp \left( \int_{t_{x,v}}^0 a(\gamma_{x,v}(s))ds \right), \quad a \in C(M),\end{equation}
where $t_{x,v}$ is the maximal backward time for the geodesic $\gamma_{x,v}$ (possibly infinite). The attenuated ray transform is the mathematical basis for the medical imaging method SPECT \cite{F86,N02b,F03}. Other applications of matrix weighted ray transforms are described in the introduction of the article \ref{paper2}. More details and references can be found from \cite{IM18}.

We list some positive injectivity results next: \begin{itemize}
\item If $(M,g)$ is a compact Riemannian manifold of $\dim(M) \geq 3$ with a strictly convex boundary and admits a smooth strictly convex function, and $W \in C^\infty(SM; GL(k,\C))$, then $I_W$ is injective for smooth functions \cite{PSUZ16}.

\item Let $(M,g)$ be a simple manifold of $\dim(M) = 2$. Let $a \in C^\infty(M)$ be a complex function and $I^a = I_{w_a}$ the attenuated ray transform with the weight $w_a$. Suppose that $F(x,v) = f(x) + \alpha(x,v)$ is the sum of a function $f \in C^\infty(M)$ and a $1$-form $\alpha \in C^\infty(T^1M)$. If $I^a F = 0$, then $F(x,v) = ap(x)+\nabla p(x,v)$ for some $p \in C^\infty(M)$ with $p|_{\partial M}=0$ \cite{SU11}. The result generalizes to the matrix weighted case where the matrix weight is the sum of a smooth unitary connection and a smooth skew-Hermitian matrix function \cite{PSU12}, and to higher dimensions if $K\leq 0$ \cite{GPSU16,PS18}.

\item If $M$ has a strictly convex boundary and $w \in C(SM)$, then $I_wf$ determines the boundary jet of a smooth function \cite{I14}. Hence, $I_w$ is injective for analytic functions. This result is based on a local argument and generalizes to the matrix weighted case straightforwardly even though it is not stated in \cite{I14}.

\item Many positive results are known in Euclidean spaces. If $n \geq 2$, $w$ is smooth, and has a rotation invariance \cite{Q83} or $w$ is real analytic \cite{BQ87}, then $I_w$ is injective. If $n \geq 3$ and the weight is regular enough ($C^{1,\alpha}$ is sufficient for example), then $I_w$ is injective \cite{MQ85,F86,I16}.
\end{itemize}

There are two important counterexamples for uniqueness in Euclidean spaces \cite{B93, GN17}. The counterexample in \cite{B93} gives a construction of a smooth weight so that the kernel of $I_w$ is nontrivial on the unit disk of the plane. The counterexample in \cite{GN17} gives a construction of a $\alpha$-Hölder continuous rotation invariant weight (in the sense of \cite{Q83}) in $\R^n$, $n\geq 2$, for some small $\alpha > 0$, so that the kernel of $I_w$ is nontrivial. This also gives a counterexample to the result of \cite{Q83} if the weight is not regular enough.

In the article \ref{paper2}, we restrict our study to the class of piecewise constant functions. We show that under this assumption continuity of a matrix weight is sufficient for showing that $I_W f = 0$ implies $f=0$. This result is valid for manifolds of $\dim(M) \geq 2$ that admit a strictly convex function.

\subsection{Geodesic ray transform on closed manifolds}

Suppose that $(M,g)$ is a closed Riemannian manifold with $\dim(M) \geq 2$. Let $\Gamma_c \subset \Gamma$ be the set of closed unit speed geodesics. Let $\tau_\gamma$ be the smallest period of $\gamma \in \Gamma_c$. The \textit{geodesic ray transform on a closed manifold} is defined by \begin{equation}\label{eq:pergrt}If(\gamma) = \int_0^{\tau_\gamma} f(\gamma(t))dt.\end{equation} This definition can be generalized to the functions on $SM$ as well.

There is again a vast literature on the geodesic ray transforms of this type in general. A lot is known for flat tori, Lie groups and other symmetric spaces \cite{I15, I16Lie, H99, Hel13}. More generally, the geodesic ray transform has been studied on Anosov surfaces and manifolds of negative curvature \cite{PSU15,GPSU16,PSU14anov}. It has applications to the \textit{spectral rigidity problem} which asks if the spectrum of the Laplace-Beltrami operator determines the metric up to a natural gauge \cite{GK80a, GK80}. 

A historically interesting fact is that the geodesic ray transform of $S^2$, called the \textit{Funk transform}, was studied for the first time by Hermann Minkowski in the early 1900s \cite{M04} and by Paul Funk a few years later \cite{F13,F15}, about a decade before the first studies of Radon on $\R^2$. The injectivity result on $S^2$ states that a symmetric function can be uniquely determined from its line integrals over great circles \cite{F15}.

In the article \ref{paper3}, we study the ray transform of closed geodesics in the special case of the flat torus $(\T^2,g_E)$. Our arguments in \ref{paper3} are specialized to the case of the flat tori and based on rather simple analysis of Fourier series. The work \ref{paper3} has applications in computational reconstructions from practical X-ray data since the geometry is flat. These results are further generalized to the periodic $d$-plane Radon transforms on $(\T^n,g_E)$ in the article \ref{paper4}. These generalizations require suitable weighted Sobolev spaces on the image side, and give another view of the theorems in \ref{paper3} in terms of weighted Sobolev spaces.

\subsection{Related open problems}\label{sec:openproblems}

We list here some important open problems in integral geometry \cite{PSU14,IM18}:

\begin{enumerate}[i)]
\item Is $I_m$ s-injective for $m \geq 2$ if $(M,g)$ is a simple manifold and $\dim(M) \geq 3$? \label{sinj3}
\item Is $I_m$ s-injective for $m \geq 0$ if $(M,g)$ is a nontrapping manifold and $\dim(M) \geq 2$? \label{nontrap}
\item If $(M,g)$ is a simple or a nontrapping manifold and $\dim(M) \geq 3$, does there exists a strictly convex function? \label{strict}
\item Is $I_m$ s-injective for $m \geq 0$ if $(M,g)$ has a strictly convex boundary and a strictly convex function, and $\dim(M) = 2$? \label{sinj2}
\item Is the attenuated geodesic ray transform $I^a$ injective if $(M,g)$ is a simple manifold and $\dim(M) \geq 3$? \label{atte}
\item Is the attenuated geodesic ray transform $I^a$ injective if $(M,g)$ is a nontrapping manifold and $\dim(M) \geq 2$?
\item Is the class of simple metrics of $M$ with $\dim(M) \geq 3$ \textit{boundary distance rigid}: Suppose that $g$ and $h$ are simple metrics on $M$. Does $d_g|_{\partial M \times \partial M} = d_h|_{\partial M \times \partial M}$ imply that $g = \phi_* h$ for some diffeomorphism $\phi: M \to M$ with $\phi|_{\partial M} = \text{Id}$? \label{brp}
\end{enumerate}

If one can solve one of the corresponding problems for nontrapping manifolds with a positive answer, then this would solve the corresponding problem for simple manifolds. Vice versa, counterexamples for simple manifolds would serve as counterexamples for nontrapping manifolds. The positive answer to the question \ref{nontrap} was conjectured in \cite{PSU13} when $\dim(M) = 2$, and the problem \ref{sinj2} is equivalent to \ref{nontrap} in this case \cite{PSUZ16}. A positive answer to \ref{strict} in the case of simple manifolds would imply a positive answer to the boundary rigidity problem \ref{brp} \cite{UV16,SUV17} and the injectivity problem \ref{atte} \cite{PSUZ16}. The positive answer to the problem \ref{brp} was conjectured by Michel in 1981 \cite{M81}, and was proved when $\dim(M) = 2$ by Pestov and Uhlmann in 2005 \cite{PU05}. As far as the author knows, there do not exist positive theorems or counterexamples to the precise statements of the problems in this list.

Injectivity of the geodesic ray transform with a smooth weight is also open on simple manifolds of $\dim(M) \geq 3$. In $\dim(M) = 2$, a positive answer cannot be obtained due to the smooth counterexample of Boman \cite{B93} on Euclidean plane. Minimal regularity assumptions of the weights for which injectivity of $I_w$ holds is also an open question in $\R^n$, $n \geq 3$ \cite{GN17,I16}. For example, is $I_w$ injective on smooth functions of the closed unit ball of $\R^n$, $n \geq 3$, if $w$ is Lipschitz continuous?

If a Riemannian manifold $M$ is assumed to be noncompact, then there are many results in symmetric geometries, but several questions of integral geometry are yet unstudied in more general geometries. The article \ref{paper1} and the work \cite{GGCS17} contain the only s-injectivity results, that the author is aware of, when special symmetries such as a constant curvature is not assumed. A further discussion on the geodesic ray transform on noncompact manifolds is given in section \ref{paper1intro} of the thesis.

\section{Main results}\label{sec:main}

\subsection{S-injectivity of the geodesic ray transform on Cartan-Hadamard manifolds, \ref{paper1}}\label{paper1intro}

We begin by introducing some notations and definitions. We then state our main results in the article \ref{paper1} and discuss earlier works in tensor tomography on noncompact manifolds. We finish this section by giving an outline of the used methods and arguments.

Let $(M,g)$ be a Cartan-Hadamard manifold. Fix a point $o \in M$. If $\eta > 0$ and $f \in C(M)$, we say that $f$ \textit{decays exponentially} and denote that $f \in E_\eta(M)$ if \begin{equation}\abs{f(x)} \leq Ce^{-\eta d(x,o)}\quad \text{for some $C > 0$},\end{equation} and $f$ \textit{decays polynomially} and denote that $f \in P_\eta(M)$ if \begin{equation}\abs{f(x)} \leq C(1+d(x,o))^{-\eta}\quad  \text{for some $C > 0$}.\end{equation}
Let $f \in C^1(M)$. We denote $f \in E_\eta^1(M)$ if $\abs{f(x)} + \abs{\nabla f(x)} \in E_\eta(M)$, and $f \in P_\eta^1(M)$ if $\abs{f(x)} \in P_\eta(M)$ and $\abs{\nabla f(x)} \in P_{\eta+1}(M)$.

Let $f,h \in C^1(T^mM)$. We define the standard inner product for $m$-tensors on $T_xM$ by \begin{equation}g_x(f,h) := g^{j_1k_1}(x)\cdots g^{j_mk_m}(x)f_{j_1\cdots j_m}(x)h_{k_1 \cdots k_m}(x).\end{equation} The norm is defined by $\abs{f}_g := \sqrt{g(f,f)}$ and defines a mapping $M \to \R$. If $f \in C^1(S^mM)$, then we write $f \in E_\eta(M)$ if $\abs{f}_g \in E_\eta(M)$, and $f \in E_\eta^1(M)$ if $\abs{f}_g \in E_\eta(M)$ and $\abs{\nabla f}_g \in E_\eta(M)$. We define analogously the sets $P_\eta(M)$ and $P_\eta^1(M)$ for tensors.

In [\ref{paper1}, Lemma 4.1], we show that $I_mf$ is well defined if $f \in P_\eta$ for some $\eta > 1$. Since $M$ is noncompact and every geodesic has infinite length, this must be shown. It is also straightforward to argue that the kernel of $I_m$ contains symmetric tensors of the form $\sigma_m(\nabla f)$ such that $f \in C(S^{m-1}M)$ and $f$ has suitable decay at infinity. We are ready to state our main results on s-injectivity of $I_m$ on Cartan-Hadamard manifolds.

\begin{theorem}[\ref{paper1}, Theorem 1.1]\label{thm:A1} Let $(M,g)$ be a Cartan-Hadamard manifold of dimension $n \geq 2$ with $-K_0 \leq K \leq 0$ for some $K_0>0$. Let $f \in E_\eta^1(M)$ be a symmetric $m$-tensor field for some $\eta > \frac{n+1}{2}\sqrt{K_0}$. If $I_mf =0$, then $f = \sigma_m(\nabla h)$ for some symmetric $(m-1)$-tensor field $h$ such that $h\in E_{\eta-\epsilon}(M)$ for any $\epsilon > 0$. (If $m = 0$, then $f \equiv 0$.)
\end{theorem}

\begin{theorem}[\ref{paper1}, Theorem 1.2]\label{thm:A2} Let $(M,g)$ be a Cartan-Hadamard manifold of dimension $n \geq 2$ and assume that $\mathcal{K} \in P_\kappa(M)$ for some $\kappa > 2$. Let $f \in P_\eta^1(M)$ be a symmetric $m$-tensor field for some $\eta > \frac{n+2}{2}$. If $I_mf =0$, then $f = \sigma_m(\nabla h)$ for some symmetric $(m-1)$-tensor field $h$ such that $h\in P_{\eta-1}(M)$. (If $m = 0$, then $f \equiv 0$.)
\end{theorem}

These theorems extend the earlier results in \cite{Leh16} where the same problem was studied in the case of functions ($m=0$) and $n = 2$. We remark that the proof of \cite[Lemma 4.6]{Leh16} is incomplete, and hence, the theorems cannot be used as stated in \cite{Leh16}. Theorems \ref{thm:A1} and \ref{thm:A2} here are proved by a different method and thus the corresponding lemma is not required. However, there might be a possibility to find better lower bounds for $\eta$ in theorems \ref{thm:A1} and \ref{thm:A2} by combining arguments of \cite{Leh16} and \ref{paper1} carefully.

The geodesic ray transform for functions on noncompact manifolds has been studied before in Euclidean and hyperbolic spaces \cite{Hel94,H99,Jen04}, and for vector fields in \cite{B05}. In these works, the regularity and decay conditions are sharper than those in theorems \ref{thm:A1} and \ref{thm:A2}. Differentiability is not needed but similar decay conditions for the function itself is required with slightly better lower bounds for $\eta$. There exist counterexamples if one does not assume a decay condition \cite{Zal82,AG93}. Theorem \ref{thm:A1} resembles the hyperbolic results and theorem \ref{thm:A2} the Euclidean. Our differentiability assumption comes from the method of proof that is based on the Pestov identity.

There are also works in noncompact spaces of constant curvature and noncompact homogeneous spaces \cite{H99,Hel13}. Theorems \ref{thm:A1} and \ref{thm:A2} are the first results on the geodesic ray transform of noncompact manifolds without special symmetries, which the author is aware of.

There is a recent related work \cite{GGCS17} where s-injectivity for $I_m$ for $m = 0,1$ was shown in the case of asymptotically hyperbolic manifolds without conjugate points and with hyperbolic trapped set. It was also shown there that if additionally $K \leq 0$, then $I_m$ is s-injective for any $m \geq 0$. The results in \ref{paper1} are not included in \cite{GGCS17}, and vice versa. Of course, there are geometries which satisfy the assumptions of the both works.

\begin{proof}[Outline of the proof of theorems \ref{thm:A1} and \ref{thm:A2}]

Let $f \in P_\eta(M)$ be a symmetric $m$-tensor field. One defines the function
\begin{equation}u^f(x,v) := \int_0^\infty \lambda_mf(\gamma_{x,v}(t),\dot{\gamma}_{x,v}(t))dt.\end{equation} A simple calculation shows that \begin{equation}u^f(x,v) + (-1)^m u^f(x,-v) = I_mf(x,v) < \infty.\end{equation} We write here $f = \lambda_m f$ to keep notation shorter. It can be calculated that $Xu^f = -f$ where $X$ is the geodesic vector field. Now one needs to understand the system $Xu = -f$ when $f$ is a symmetric $m$-tensor such that $I_mf = 0$ and $f$ satisfies the assumptions of theorem \ref{thm:A1} or theorem \ref{thm:A2}.

We list the main ideas next:
\begin{enumerate}[i)]
\item The goal is to show that $f = -Xu^f = \lambda_m(\sigma_m \nabla U)$ for some $U \in C(S^{m-1}M)$ with right decay properties. \label{proof1}
\item If $M$ is a compact manifold with boundary and $K \leq 0$, then the Pestov identity can be used to show \ref{proof1}. This follows from a contraction property of the Beurling transform on manifolds of nonpositive sectional curvature \cite{PSU15}. \label{proof2}
\item The energy estimates of the step \ref{proof2} in compact manifolds that involve only terms up to the first order derivatives can be extended to $H^1(SM)$ when $M$ is a complete manifold with $K \leq 0$. These $H^1(SM)$ extensions of the energy estimates and the final argument to show \ref{proof1} are done in [\ref{paper1}, Section 5].
\item Hence, we need to show that $u^f \in H^1(SM)$ under the assumptions of theorems \ref{thm:A1} and \ref{thm:A2}. The core part of this is done in [\ref{paper1}, Section 4] by estimating growths of Jacobi fields on Cartan-Hadamard manifolds.
\end{enumerate}

We next explain some of the details. Showing that $u^f \in H^1(SM)$ is a bit tricky and our argument uses geometric estimates for growths of Jacobi fields and the decay assumptions of $f$. The idea could be summarized as follows: the faster the geodesics spread the faster the functions (and derivatives) should decay to make $L^2$ estimates work because of the growth rate of volumes of balls (cf. [\ref{paper1}, Lemmas 4.8 and 5.4]).

One can orthogonally split the gradient of $SM$ as \begin{equation}\nabla_{SM}u = (Xu)X + \hnabla u + \vnabla u \label{eq:nabla}\end{equation} where $X$ and $\hnabla$ represents \textit{horizontal derivatives} with respect to $x$ and $\vnabla$ \textit{vertical derivatives} with respect to $v$. These and other geometric preliminaries are given in [\ref{paper1}, Section 3]. The most technical part is the proof of [\ref{paper1}, Lemma 4.7]. In that lemma, we first show that $u^f$ is locally Lipschitz and then estimate the components (\ref{eq:nabla}) of the gradient $\nabla_{SM} u^f$ for a.e. $(x,v) \in SM$ based on our Jacobi field estimates. This implies that $u^f \in H^1(SM)$ [\ref{paper1}, Lemma 5.4].

The rest of the argument uses estimates and methods developed in \cite{PSU15}. Details of the spherical harmonics decomposition of $L^2(SM)$ are given in \cite{GK80,DS11}. Let $H_k(SM)$ be the eigenspace for the eigenvalue $k(k+n-2)$ of the spherical Laplacian. One can split the geodesic vector field $X = X_+ + X_-$ into two parts so that $X_+: \Omega_k \to H_{k+1}(SM)$ and $X_-: \Omega_k \to H_{k-1}(SM)$ where $\Omega_k = H_k(SM) \cap H^1(SM)$. We can show this by proving the estimate
\begin{equation}\norm{X_+u}^2 + \norm{X_-u}^2 \leq \norm{Xu}^2 + \norm{\hnabla u}^2\end{equation} for $u \in H^1(SM)$ [\ref{paper1}, Lemma 5.1]. This part of the proof requires the Pestov identity and estimates based on the contraction property of the Beurling transformation from \cite{PSU15}. If $u \in H^1(SM)$, it follows that the spherical harmonics decomposition has the form \begin{equation} u = \sum_{k=0}^\infty u_k, u_k \in \Omega_k,\end{equation} where the series converges in $L^2(SM)$. We can now conclude that if $u \in H^1(SM)$, then $\norm{X_+ u_k} \to 0$ as $k \to \infty$ [\ref{paper1}, Corollary 5.2].

Since symmetric $m$-tensors have only terms up to degree $m$ in their spherical harmonic decomposition, we get
\begin{equation}-\sum_{k=0}^m f_k = -f = Xu^f = X_+u^f + X_-u^f.\label{eq:joku}\end{equation} The rest of the proof follows from the formula (\ref{eq:joku}) and [\ref{paper1}, Corollary 5.2 and Lemma 5.3] by following arguments from \cite{PSU13,PSU15}. The final step is to straightforwardly estimate decay of the elements of the kernel. These details are given in [\ref{paper1}, Proof of theorems 1.1 and 1.2.].\end{proof}




\subsection{On the geodesic ray transform with matrix weights for piecewise constant functions, \ref{paper2}}

The geodesic ray transform for piecewise constant functions was studied on the manifolds that admit a strictly convex function in \cite{ILS17}. The work \cite{ILS17} was motivated by the fact that injectivity of the geodesic ray transform is an open problem for nontrapping manifolds. If $n = 2$, then a manifold with strictly convex boundary is nontrapping if and only if it has a strictly convex function (see \cite[Section 2]{PSUZ16} for details and references). The main result of \cite{ILS17} was to show that $If = 0$ implies $f\equiv0$ if $f$ is a piecewise constant function on $M$. Reconstruction of a piecewise constant function from $If$ was studied recently in \cite{Leb19}.

Piecewise constant functions are defined according to the definition of \cite{ILS17}. We recall this definition next. A \textit{regular tiling} of a manifold is a collection of regular $n$-simplices which cover the manifold, whose interiors are disjoint, and whose boundaries intersect nicely \cite[Section 2.1]{ILS17}. A function $f: M \to \C^k$ is called \textit{piecewise constant} if there exists a regular tiling $\{\Delta_1,\dots,\Delta_N\}$ such that $f|_{\text{Int}(\Delta_i)}$ is constant for any $i \in \{1,\dots,N\}$ and $f \equiv 0$ elsewhere.

The main result of the article \ref{paper2} generalizes the main result of \cite{ILS17} to the matrix weighted case, analogous to the problem studied in \cite{PSUZ16} for smooth functions and weights in dimensions $n \geq 3$. We denote by $\text{Mon}(\C^k,\C^m)$ the space of injective linear maps $\C^k \to \C^m$.

\begin{theorem}[\ref{paper2}, Theorem 1.1]\label{thm:thmpaper2} Let $(M,g)$ be a compact nontrapping Riemannian manifold with strictly convex smooth boundary and $W \in C(SM;\text{Mon}(\C^k,\C^m))$. Let either \begin{enumerate}[(a)]
\item $\dim(M) = 2$, or
\item $\dim(M) \geq 3$ and $(M,g)$ admits a smooth strictly convex function.
\end{enumerate} If $f: M \to \C^k$ is a piecewise constant vector-valued function and $I_Wf = 0$, then $f \equiv 0$.
\end{theorem}

\begin{remark} Piecewise constant functions do not form a vector space under the definition used in the study \cite[Remark 2.7]{ILS17}. Hence, injectivity follows only if the tiling of the piecewise constant function are known beforehand. It is an open problem how to determine the tiling of a piecewise constant function from the data $If$.
\end{remark}

The proof of theorem \ref{thm:thmpaper2} is strongly based on the method developed in \cite{ILS17}. We show that locally the matrix weighted geodesic ray transform data can be reduced to the data of the geodesic ray transform without weight [\ref{paper2}, Lemma 2.4 and Lemma 2.5]. We remark that this reduction does not work for general functions but it works for piecewise constant functions. Local injectivity of the geodesic ray transform for piecewise constant functions was shown in \cite{ILS17}. The layer stripping argument of \cite{ILS17}, using a strictly convex function, allows one to go from the local uniqueness result to the global uniqueness result [\ref{paper2}, Theorem 2.6].

\subsection{Theory of Tikhonov regularized reconstructions from the X-ray transform data on the flat 2-torus, \ref{paper3}}

The geodesic ray transform on the flat torus $\T^2 := \R^2 / \Z^2$ is defined for the closed geodesics. A geodesic is closed on $\T^2$ if and only if its directional vector is a multiple of an integer vector. Instead of unit-speed parametrization of geodesics, we parametrize geodesics so that each closed geodesic has the period $1$. This is convenient since the $1$-periodic geodesics are of the form \begin{equation}\gamma_{x,v}(t) = \pi(x+tv),\quad x \in \R^2, v\in \Z^2 \setminus 0, t \in [0,1]\end{equation}
where $\pi: \R^2 \to \T^2$ is the quotient mapping. Clearly, if $\pi(x) = \pi(y)$, then $\gamma_{x,v} = \gamma_{y,v}$ for any $v\in\Z^2\setminus 0$.

Hence, the (geodesic) X-ray transform on $\T^2$ can be defined by
\begin{equation}If(x,v) = \int_0^1 f(\gamma_{x,v}(t))dt\end{equation}
for continuous functions. We remark that this definition actually scales the data (\ref{eq:pergrt}) by the factor $\abs{v}^{-1}$. However, there is one-to-one correspondence between the both definitions of $I$ on $\T^2$. This definition extends to the dual space of smooth functions, called \textit{distributions} and denoted by $\mathcal{T}'$, since $If(\cdot,v)$ is formally $L^2(\T^2)$ self-adjoint for every fixed $v \in \Z^2\setminus 0$. For further details see \cite{I15, I17} or [\ref{paper3}, Section 2.1].

Injectivity of $I$ on tori is well understood and it has been studied earlier in \cite{S82,AR11,A11,I15}. The main contributions of \ref{paper3} are related to reconstruction, better understanding of functional properties, and numerical simulations that demonstrate applicability of the method in CT imaging. It is described in [\ref{paper3}, Section 2.3] and \cite[Chapter 3]{I17} how practical X-ray data of a compactly supported object on $\R^2$ can be mapped into X-ray data on $\T^2$.

One has the \textit{Fourier series decomposition}
\begin{equation}f(x) = \sum_{k\in \Z^2} \hat{f}(k)e^{2\pi ik \cdot x}, \quad \hat{f}(k) := f(e^{-2\pi i k \cdot x}),k \in \Z^2,\end{equation}
for any $f \in \mathcal{T}'$. It was shown in \cite[Eq. (9)]{I15} that for any $f \in \mathcal{T}'$ the identity
\begin{equation}\label{eq:ilmavirta}\widehat{If}(k,v) = \begin{cases} \hat{f}(k) & k \cdot v = 0\\
0 & k \cdot v \neq 0\end{cases}\end{equation}
holds. This gives a reconstruction formula for $f$ from the data $If$ and shows injectivity. In the work \ref{paper3}, we have studied consequences of this formula further and implemented a reconstruction algorithm based on our new findings.

We state and describe our main theorems in \ref{paper3} next. Our first theorem simplifies the reconstruction formula (\ref{eq:ilmavirta}) for integrable functions. This simplification results better computational efficiency since the dimension of the integrals (\ref{eq:ilmavirta}) are reduced by one.

\begin{theorem}[\ref{paper3}, Theorem 1]
\label{thm:invfor}
Suppose that $f \in L^1(\T^2)$. Let $k \in \Z^2$. If $k,v \neq 0$ and $v \bot k$, then
\begin{equation}
\label{eq:rec1}
\hat{f}(k) = \begin{cases} \int_0^1 I_vf(0,y)\exp(-2 \pi ik_2y)dy,\quad k_2 \neq 0 \\
\int_0^1 I_vf(x,0)\exp(- 2 \pi i k_1x)dx,\quad k_1 \neq 0.\end{cases}\end{equation} If $k = 0$, then \begin{equation}\label{eq:rec0}\hat{f}(k) = \int_0^1 I_{(1,0)}f(0,y)dy = \int_0^1 I_{(0,1)}f(x,0)dx.
\end{equation}
\end{theorem}

This theorem can be proved by a change of coordinates and Fubini's theorem. We gave two proofs in [\ref{paper3}, Section 2.2]. The first proof gives a new proof of injectivity of $I$ on $\T^2$. The second proof uses the formula (\ref{eq:ilmavirta}) directly. A slightly more general statement is actually proved in [\ref{paper3}, Theorem 8].

Our next two theorems are about regularization. We need to first introduce a suitable Sobolev space structure on the image side. Let $Q\subset\Z^2$ be such that every nonzero $v\in\Z^2$ is an integer multiple of a unique element in~$Q$. This set can be naturally identified with the rational projective space~$\mathbb P^1$. The X-ray transform  takes a function on~$\T^2$ to a function on $\T^2\times Q$. 

\begin{remark} There is a connection between X-ray tomography with partial data and Schanuel's theorem \cite{S64} on heights of projective spaces [\ref{paper3}, Section 2.6.2]. In particular, the number of directions $v \in \Z^2 \setminus 0$ needed in the reconstruction of the Fourier coefficients of $f$ in $B_{\ell^\infty}(0,R)$ from $I_vf$ can be estimated using Schanuel's theorem.
\end{remark}

We use the standard Sobolev scale of spaces~$H^s(\T^2)$ with the norms
\begin{equation}
\aabs{f}_{H^s(\T^2)}^2
=
\sum_{k\in\Z^2}
\vev{k}^{2s}\abs{\hat f(k)}^2,
\end{equation}
where $\vev{k}=(1+\abs{k}^2)^{1/2}$ as usual.
On $\T^2\times Q$, we define the spaces $H^{s}(\T^2\times Q)$ to be the set of functions $g\colon\T^2\times Q\to\C$ for which
\begin{enumerate}[(i)]
\item  $g(\cdot,v)\in H^s(\T^2)$ for every $v\in Q$,
\item the average of every $g(\cdot,v)$ over $\T^2$ is the same, and
\item the norm
\begin{equation}
\aabs{g}_{H^{s}(\T^2\times Q)}^2
=
\abs{\hat g(0,0)}^2+
\sum_{k\in\Z^2\setminus0}
\sum_{v\in Q}
\vev{k}^{2s}\abs{\hat g(k,v)}^2
\end{equation}
is finite.
We set $v=0$ for the Fourier term $k=0$ to emphasize that it is the same for every $v\in Q$.
We remind the reader that $0\notin Q$.
\end{enumerate}

Now, we can consider a \textit{Tikhonov minimization problem}: given some data $g\in H^r(\T^2\times Q)$, find
\begin{equation}
\label{eq:reg-min}
\argmin_{f \in H^r(\T^2)}
\left(\aabs{I f-g}_{H^r(\T^2\times Q)}^2 + \alpha\aabs{f}^2_{H^s(\T^2)}\right)
.\end{equation}
Let us define the post-processing operator~$P^s_\alpha$ to be the Fourier multiplier $(1+\alpha\vev{k}^{2s})^{-1}$ and denote by~$I^*$ the adjoint of~$I$. Formulas that define the adjoint and normal operators are proved in [\ref{paper3}, Proposition 11]. In fact, the X-ray transform is unitary as a mapping $H^s(\T^2) \to H^s(\T^2 \times Q)$ for any $s \in \R$. 

\begin{theorem}[\ref{paper3}, Theorem 2]
\label{thm:regularization}
Let $r\in\R$, $s\geq r$, and $\alpha>0$.
Suppose $g\in H^r(\T^2\times Q)$.
The unique minimizer~$f$ of the minimization problem~\eqref{eq:reg-min} corresponding to Tikhonov regularization is $f=P^{s-r}_\alpha I^*g\in H^{2s-r}(\T^2) \subset H^s(\T^2)$.
\end{theorem}

\begin{theorem}[\ref{paper3}, Theorem 3]
\label{thm:strategy}
Suppose $r,t,s,\delta\in\R$ are such that $2s+t\geq r$, $\delta\geq0$, and $s>0$.
We assume that $f\in H^{r+\delta}(\T^2)$ and $g\in H^t(\T^2\times Q)$.

Then our regularized reconstruction operator~$P^s_\alpha I^*$ gives a regularization strategy in the sense that
\begin{equation}
\lim_{\eps\to0}
\sup_{\aabs{g}_{H^t(\T^2\times Q)}\leq\eps}
\aabs{P^s_{\alpha(\eps)}I^*(I f+g)-f}_{H^r(\T^2)}
=
0,
\end{equation}
where $\alpha(\eps)=\sqrt{\eps}$.

Moreover, if $\aabs{g}_{H^t(\T^2\times Q)}\leq\eps$, $0<\delta<2s$ and $0<\alpha\leq 2s/\delta-1$, we have
\begin{equation}
\label{eq:strategy-estimate}
\aabs{P^s_\alpha I^*(I f+g)-f}_{H^r(\T^2)}
\leq
\alpha^{\delta/2s}
C(\delta/2s)
\aabs{f}_{H^{r+\delta}(\T^2)}
+
\frac{\eps}{\alpha}
,
\end{equation}
where $C(x)=x(x^{-1}-1)^{1-x}$.
\end{theorem}

A simple calculation shows that the optimal rate of convergence is obtained if the regularization parameter is chosen so that $\alpha = \epsilon^\lambda$ where $\lambda = (1+\delta/2s)^{-1}$.

The proofs of the theorems are based on quite straightforward computations on the Fourier side and the formula (\ref{eq:ilmavirta}). It seems that the key theoretical finding in \ref{paper3} was the right structure on the image side. It is quite easy to see that $I$ is non-surjective between the Sobolev spaces $H^s(\T^2)$ and $H^s(\T^2 \times Q)$. Hence, the choices made for the image side Sobolev norms do not fully trivialize the problem and, instead of that, those choices describe the behavior of $I|_{H^s(\T^2)}$.

Numerical implementation, simulations and conclusions are described in [\ref{paper3}, Sections 3--5]. A short discussion of typical numerical methods in CT imaging is given in [\ref{paper3}, Section 1.2]. We do not repeat the details or discussions here. We did not perform tests with measured X-ray laboratory data. This would be the next step towards practical CT imaging based on the reconstruction method on the flat torus.

\subsection{Fourier analysis of periodic Radon transforms, \ref{paper4}}

The article \ref{paper4} studies the periodic $d$-plane Radon transforms on $\T^n := \R^n /\Z^n$ when $1 \leq d \leq n-1$ and $n \geq 2$. If $n =2$ and $d=1$, then the $d$-plane Radon transform is the X-ray transform studied in the article \ref{paper3}. The periodic Radon transforms have been applied in other mathematical tomography problems earlier: the broken ray transform on boxes \cite{I15}, the geodesic ray transform on Lie groups \cite{I16Lie}, tensor tomography on periodic slabs \cite{IU18}, and the ray transforms on Minkowski tori \cite{I18}.

We generalize the main theorems in \ref{paper3} into higher dimensions [\ref{paper4}, Theorems 1.4 and 1.5, Proposition 3.1]. We do not restate these statements here. We state here results on the adjoint and normal operators and the stability estimates. We also introduce a new inversion formula which might be of a practical interest due to its simplicity.

We begin by introducing necessary mathematical preliminaries. Suppose that $f \in \mathcal{T} := C^\infty(\T^n)$, then we define the \textit{$d$-plane Radon transform} of $f$ by
\begin{equation}
R_df(x,A) := \int_{[0,1]^d}f(x+t_1v_1+\cdots+ t_dv_d)dt_1\dots dt_d
\end{equation}
where $A = \{v_1,\dots,v_d\}$ is a set of $d$ linearly independent integer vectors $v_i \in \Z^n$. 

It can be shown that $A$ spans a periodic $d$-plane on $\T^n$. On the other hand, if $A$ and $B$ span the same periodic $d$-plane on $\T^n$, then $R_df(x,A) = R_df(x,B)$ for any $x \in \T^n$. Let $\Gr(d,n)$ denote the collection of $d$-dimensional subspaces of $\Q^n$. These spaces are called \textit{Grassmannians}. For any element in $\Gr(d,n)$ there exists a basis of integer vectors. Hence, we may define $R_df: \Gr(d,n) \to \mathcal{T}$ using bases of integer vectors as representatives of elements in $\Gr(d,n)$. The definition of $R_d$ extends to the periodic distributions $\mathcal{T}'$ using the duality and the fact that $R_d(\cdot,A): \mathcal{T} \to \mathcal{T}$ is formally $L^2$ self-adjoint for any fixed $A \in \Gr(d,n)$. Let us denote $R_{d,A}f = R_df(\cdot,A)$ for any $f \in \mathcal{T}'$.

Next, we define suitable structures for the data spaces such that the images of the Bessel potential spaces $L_s^p(\T^n)$ under $R_d$ are contained into the data spaces. Let $p,l \in [1,\infty]$ and $s \in \R$. We define the \textit{Bessel potential norms} as
\begin{equation}\begin{split}\norm{f}_{L_s^p(\T^n)} &= \norm{\sum_{k \in \Z^n} \vev{k}^s \hat{f}(k)e^{2\pi i k \cdot x}}_{L^p(\T^n)}, \\
\norm{f}_{H^s(\T^n)} &= \sqrt{\sum_{k \in \Z^k} \vev{k}^{2s}\abs{\hat{f}(k)}^2}\end{split}\end{equation} where $\vev{k}=(1+\abs{k}^2)^{1/2}$ as usual. The space $L_s^p(\T^n) \subset \mathcal{T}'$ consists of all $f \in \mathcal{T}'$ with $\norm{f}_{L_s^p(\T^n)} < \infty$. If $p = 2$, then $H^s(\T^n) = L_s^p(\T^n)$. One has equivalently that $f \in L_s^p(\T^n)$ if and only if $(1-\Delta)^{s/2}f \in L^p(\T^n)$ and $f\in \mathcal{T}'$.

Let us denote $X_{d,n} := \T^n \times \Gr(d,n)$ to keep our notation shorter. Let $w: \Z^n \times \Gr(d,n) \to (0,\infty)$ be a weight function such that $w(\cdot,A)$ is at most of polynomial decay for any fixed $A \in \Gr(d,n)$ (see [\ref{paper4}, Section 2.2] for the definition). We say that a function $g: X_{d,n} \to \C$ belongs to $L_{s}^{p,l}(X_{d,n}; w)$ with $1 \leq l < \infty$ if the norm
\begin{equation} \norm{g}_{L_{s}^{p,l}(X_{d,n}; w)}^l := \sum_{A \in \Gr(d,n)} \norm{g(\cdot,A)}_{L_s^p(\T^n; w(\cdot,A))}^l
\end{equation}
is finite and $g(\cdot,A) \in \mathcal{T}'$ when $A \in \Gr(d,n)$. Similarly, if $l = \infty$, we define
\begin{equation} \norm{g}_{L_{s}^{p,\infty}(X_{d,n}; w)} := \sup_{A \in \Gr(d,n)} \norm{g(\cdot,A)}_{L_s^p(\T^n; w(\cdot,A))}
\end{equation}
If $p,l = 2$, then the norm is generated by the corresponding inner product. The spaces $L_s^{p,l}(X_{d,n};w)$ are Banach spaces [\ref{paper4}, Lemma 2.1].

We have introduced weighted structures since most of the theorems in \ref{paper4} would have been unreachable without such structures when $d < n-1$. If $d=n-1$, then the analysis of \ref{paper3} using slightly different data spaces generalizes nicely without weights. It is explained in the article \ref{paper4} how the results in \ref{paper3} can be obtained from the results in \ref{paper4}. We construct weights that satisfy the assumptions of our theorems in [\ref{paper4}, Section 2.3].

We state some of the main results in \ref{paper4} next.

\begin{theorem}[\ref{paper4}, Theorem 1.1]\label{thm:adjoint} Let $s \in \R$ and suppose that there exists $C_w > 0$ such that \begin{equation}\sum_{A \in \Omega_k} w(k,A)^2 \leq C_w^2, \quad \Omega_k := \{\, A \in \Gr(d,n) \,;\, k\bot A\,\}\end{equation}
for any $k \in \Z^n$. Then the adjoint of $R_{d}: H^s(\T^n) \to L_s^{2,2}(X_{d,n};w)$ is given by
\begin{equation}\widehat{R_d^*g}(k) = \sum_{A \in \Omega_k} w(k,A)^2\hat{g}(k,A)\end{equation}
and the normal operator $R_d^*R_d: H^s(\T^n) \to H^s(\T^n)$ is the Fourier multiplier $W_k := \sum_{A \in \Omega_k} w(k,A)^2$. In particular, the mapping  $F_{W_k^{-1}}R_d^*: R_d(\mathcal{T}') \to \mathcal{T}'$ is the inverse of $R_d$.
\end{theorem}

Theorem \ref{thm:adjoint} generalizes [\ref{paper3}, Proposition 11] into higher dimensions and implies the following results on stability.

\begin{corollary}[\ref{paper4}, Corollary 1.2]\label{cor:h2results}Suppose that the assumptions of theorem \ref{thm:adjoint} hold, and that there exists $c_w > 0$ such that $W_k \geq c_w^2$ for any $k \in \Z^n$. \begin{enumerate}[label=(\roman*)]
\item  Then $F_{W_k^{-1}}R_d^*: L_s^{2,2}(X_{d,n};w) \to H^s(\T^n)$ is $1/c_w$-Lipschitz.\label{item:h1prop1}
\item Let $f \in \mathcal{T}'$. Then
\begin{equation}\label{eq:normalopcor}\norm{f}_{H^s(\T^n)} \leq \frac{1}{c_w}\norm{R_df}_{L_s^{2,2}(X_{d,n};w)}.\end{equation}\label{item:h1prop2}
\item Let $\tilde{w}(k,A) = \frac{w(k,A)}{\sqrt{W_k}}$ and $p \in [1,\infty]$. Then $R_d^{*,\tilde{w}}R_df = f$ and $\norm{f}_{L_s^p(\T^n)} = \norm{R_d^{*,\tilde{w}}R_df}_{L_s^p(\T^n)}$ for any $f \in \mathcal{T}'$.\label{item:Lspstab}
\end{enumerate}
\end{corollary}

Other stability estimates on $L_s^p(\T^n)$ are given in terms of $R_df$ in [\ref{paper4}, Proposition 4.3]. Those stability estimates follow from corollary \ref{cor:h2results} and the Sobolev inequality on $\T^n$. This method requires additional smoothness of $R_df$ in order to control the norm of $f$ due to the use of the Sobolev inequality. The stability estimates in \ref{paper4} are new in any dimension, and different than the stability estimates in \cite{I15}.

\begin{theorem}[\ref{paper4}, Theorem 1.3]\label{thm:recwithoutFourier} Suppose that $f \in \mathcal{T}'$. Let $w: \Z^n \times \Gr(d,n) \to \R$ be a weight so that
\begin{equation}\sum_{A \in \Omega_k} w(k,A) = 1, \quad \Omega_k := \{\, A \in \Gr(d,n) \,;\, k\bot A\,\}\end{equation} and the series is absolutely converging for any $k \in \Z^n$ (the weight does not have to generate a norm or have at most of polynomial decay). Then
\begin{equation}(f,h) = \sum_{A \in \Gr(d,n)}(F_{w(\cdot,A)}R_{d,A}f,h), \quad \forall h \in \mathcal{T}.\end{equation}
Moreover, 
if $f$ has zero average and $d = n-1$, then \begin{equation}\label{eq:recwithoutFourier0mean}f = \sum_{A \in \Gr(d,n)} R_{d,A}f.\end{equation}
\end{theorem}

Theorem \ref{thm:recwithoutFourier} gives a new reconstructive formula for the inverse of $R_d$. The case $d = n-1$ is especially interesting since it does not involve any filtering, and averages are simple to reconstruct and filter out from $R_df$. The proof of theorem \ref{thm:recwithoutFourier} follows easily from the higher dimensional version of the formula (\ref{eq:ilmavirta}) proved in \cite{I15}.

\bibliographystyle{abbrv}
\bibliography{sample}


\end{document}